\documentclass[12pt]{article}
\usepackage{hyperref}
\usepackage{amsmath}
\usepackage{ntheorem}
\usepackage{mathtools}
\usepackage{epsf}
\usepackage{theorem}
\usepackage{indentfirst}
\usepackage{latexsym}
\usepackage{amssymb}
\usepackage[dvips]{graphicx}
\usepackage{flafter}
\usepackage[numbers,sort&compress]{natbib}
\usepackage{caption}
\usepackage{textcomp}
\usepackage{supertabular}
\usepackage{longtable, booktabs}
\usepackage{hhline}
\usepackage{bigstrut,bigdelim}
\usepackage{rotating}
\usepackage{multicol}
\usepackage{multirow}
\usepackage{makeidx}
\usepackage{epsfig}
\usepackage{fancyhdr}
\newif\ifger
\gerfalse
\newtheorem{theorem}{Theorem}
\newtheorem{lemma}{Lemma}[section]

\newtheorem{definition}{Definition}[section]

\newtheorem{conjecture}{Conjecture}[section]

\newtheorem{proposition}{Proposition}[section]
\newtheorem{step}{Step}
\theoremstyle{break}
\newtheorem{example}{Example}[section]

\begin{document}
\baselineskip=19pt

\title{On primitivity and reduction for half-flag-transitive   block designs }
\author{Xiaoqin Zhan\footnote{Corresponding author: zhanxiaoqinshuai@126.com}\\
\small \it  School of Science, East China JiaoTong University, \\
\small \it  Nanchang, 330013, P.R. China}
\date{}
\maketitle

\begin{abstract}
Let $\mathcal{D} = (\mathcal{P}, \mathcal{B})$ be a $2$-$(v, k, \lambda)$ design, and let $G$ be a half-flag-transitive automorphism group of ${\cal D}$. In this  article,  we first establish three sufficient conditions for $G$ to be point-primitive: (i) $\lambda \geq (r, 2\lambda)^2$, (ii) $r > 4\lambda(k-2)$,  (iii) $(v-1,2k-2)\le2$.  Next, we  prove that for $\lambda \geq (r, 2\lambda)^2$, the group $G$ is either of affine type, almost simple type, or product type. Finally, we analyze the case where $G$ is of almost simple type and prove that if the socle of $G$ is a sporadic simple group then $G \cong \text{HS}$ and $\cal D$ is either the unique $2$-$(176, 128, 15240)$ design  or the unique $2$-$(176, 160, 19080)$ design. 

\medskip
\noindent{\bf Mathematics Subject Classification (2010):} 20B25, 05B25, 05B05

\medskip
\noindent{\bf Keywords:}  Half-flag-transitive; $2$-design; automorphism group; sporadic socle
\end{abstract}

\section{Introduction}
A $2$-$(v, k, \lambda)$ design $\mathcal{D}$ is a pair $(\mathcal{P}, \mathcal{B})$ with a set $\cal P$ of $v$ points and a set $\mathcal{B}$ of $b$ blocks such that each block is a $k$-subset of $\cal P$ and each pair of distinct points is contained in $\lambda$ blocks. The \textit{replication number} $r$ of $\mathcal{D}$ is the number of blocks incident with a given point We say $\mathcal{D}$ is \textit{nontrivial} if $2 < k < v-1$. All $2$-$(v, k, \lambda)$ designs in this paper are assumed to be nontrivial. 
An \textit{automorphism} of $\mathcal{D}$ is a permutation of the point set $\cal P$ which preserves the block set. The parameters $v, b, r, k, \lambda$ of a $2$-design satisfy the following conditions (see \cite[1.2, 1.9]{Handbook}):
 \begin{equation}\label{flagnum}
 vr=bk;
\end{equation}
 \begin{equation}\label{2flagnum}
 \lambda(v-1)=r(k-1);
\end{equation}
 \begin{equation}\label{fisherineq}
 b\ge v.
\end{equation}

The set of all automorphisms of $\mathcal{D}$ under composition of permutations forms a group, denoted by $\mathrm{Aut}(\mathcal{D})$. 
A subgroup $G$ of $\mathrm{Aut}(\mathcal{D})$ \textit{leaves invariant} a partition $\mathcal{C}$ of $\cal P$ if each element of $G$ permutes the parts of $\mathcal{C}$ setwise. A partition $\mathcal{C}$ is \textit{trivial} if either $\mathcal{C}$ consists of singleton sets, or $\mathcal{C} = \{\mathcal{P}\}$; and $G$ is \textit{point-primitive} if the only $G$-invariant partitions of $\mathcal{P}$ are the trivial ones. Otherwise $G$ is said to be \textit{point-imprimitive}. 
A \textit{flag} of $\mathcal{D}$ is a pair $(\alpha, B)$ where $\alpha \in \mathcal{P}$, $B \in \mathcal{B}$, and $B$ contains $\alpha$. Obviously, the number of flags is $vr$ (or $bk$). A subgroup $G$ of $\mathrm{Aut}(\mathcal{D})$ is said to be \textit{half-flag-transitive} if it acts transitively on the block set $\cal B$ and, for every block $B\in \cal B$, the setwise stabilizer  $G_B$ partitions the point set of $B$ into exactly two orbits of size $\frac{k}{2}$. Clearly, if $G$ is half-flag-transitive, then $G$ has two orbits on flag set of $\cal D$ with same size $\frac{bk}{2}$, which is half the number of flags.

Here, it should be noted that if $G$ has two orbits of size $\frac{bk}{2}$ on the set of flags, $G$ is not necessarily half-flag-transitive. The following example illustrates this case.
\begin{example} 
Let $G = \langle (3, 9, 7, 8)(4, 10, 5, 6), (1\,8\,2)(3\,4\,5)(6,10,7) \rangle \cong PSL(2,9) $ be a $2$-transitive group on the point set \( \mathcal{P}=\{1, 2,\ldots, 10\} \), and define
the point-block incidence structure \( \mathcal{D} =({\cal P},{\cal B}) \) with ${\cal B}=B_1^G \cup B_2^G$, where

  \[ B_1 = \{1, 2, 4,5\},  B_2 = \{1, 2, 3,7\} .\]
   It is easily known that \( G_{B_1} \cong G_{B_2} \cong S_4 \) is transitive on \( B_1 \) and \( B_2 \). Then \( \mathcal{D}\) is not a half-flag-transitive \(2\)-\((10, 4, 4)\) design, but \( G \) has two orbits with same size on flag set. Additionally, $\cal D $ is the union of two flag-transitive $2$-$(10, 4, 2)$ designs with the same point set $\cal P$.
\end{example}

As a generalization of the $2$-$(v, k, \lambda)$ design, we give the definition of the $2$-$(v, K, \lambda)$ design ({\it pairwise balanced design}), where $K$ is a set of positive integers such that $k \leq v$ for every $k \in K$.

\begin{definition}
Let $\cal P$ be a finite set with $v$ elements, called points, $\lambda$ be a positive integer, $K$ be a set of positive integers such that $ k \leq v$ for every $k \in K$. Let $\mathcal{B}$ be a finite collection of subsets of $\cal P$, called blocks. Then the incidence structure $\mathcal{D}=(\mathcal{P}, \mathcal{B})$ is called a $2$-$(v, K, \lambda)$ design if $|B| \in K$ for every $B \in \mathcal{B}$, and each pair of distinct points is contained in exactly $\lambda$ blocks.
\end{definition}


Let $\alpha \in \mathcal{P}$, and $k \in K$. Let $r^{(k)}_\alpha$ be the number of blocks having size $k$ through $\alpha$. Then for each $k \in K$,

\begin{equation}\label{gl2}
\sum_{k \in K} r^{(k)}_\alpha (k-1) =  \lambda (v-1).
\end{equation}

It is well-known that the Block's Lemma \cite{Block} states that a block-transitive automorphism group $G$
of a 2-design must be point-transitive but not necessarily point-primitive. However, when certain parameters of the design are appropriately restricted, the point-primitivity of $G$ can be deduced. For instance, Delandtsheer and Doyen in \cite{DD} proved that if the number of points $v>(\binom{k}{2}-1)^2$, then $G$ must be point-primitive.
In 1984, Camina and Gagen \cite{CG} gave the result that the block‐transitive automorphism group $G$ of a $2‐(v,k,1)$ design with $k\mid v$ must be flag‐transitive, moreover, $G$ is point‐primitive. As a generalization of this result,  Zhang and Zhou \cite{ZZ2022} study block‐transitive automorphism groups of nontrivial  $2$‐$(v,k,\lambda)$ designs with $(r,k)=1$. They prove that, for such a  design $\cal D$, if  \( G \leq \text{Aut}(D) \) is block‐transitive, then $G$ must be flag‐transitive, and furthermore, $G$ is point‐primitive of affine or almost simple type.

 In the case where $G$ is a flag-transitive automorphism group, there have been relatively more research achievements in this area. Flag-transitive $2‐(v,k,1)$ designs are point-primitive by \cite{HM}. Delandtsheer and Doyen \cite{DD} also proved that a flag-transitive, point-imprimitive $2$-design satisfies $v \leq (k-2)^2$. It follows immediately that, for a flag-transitive point-imprimitive $2$-$(v, k, \lambda) $ design, the number \( v \) of points cannot be too large relative to the block size \( k \).
  In 1965, Dembowski \cite{Demb1968} prove the following result:
\begin{proposition}\label{dem}
 Suppose that \(G\) is a flag-transitive automorphism group of a $2$-design $\cal D$. Then each of the following conditions implies that \(G\) is point-primitive:
\begin{enumerate}
\item[\rm(a)] \(\lambda > (r, \lambda)\big[(r, \lambda) - 1\big]\), in particular \((r, \lambda) = 1\) and \(\lambda = 1\);
\item[\rm(b)] \(( r - \lambda, k) = 1\);
\item[\rm(c)] \(r > \lambda(k - 3)\);

\item[\rm(d)] \((v - 1, k - 1) \leq 2\).
\end{enumerate}
\end{proposition}
Classifying flag-transitive 2-designs with the aforementioned parameter restrictions is a highly meaningful work.
However, the primary step in such classification is to employ the O'Nan-Scott theorem for reduction of primitive automorphism groups. Building upon foundational work in design theory and permutation group theory, several key results have been established regarding the classification of flag-transitive automorphism groups in 2-designs with specific parameter in Proposition \ref{dem}.

\textbf{Case (a):} 
The seminal work by Buekenhout et al.~\cite{BDD} demonstrated that any $2$-$(v,k,1)$ design $\mathcal{D}$ admitting a flag-transitive (and consequently point-primitive) automorphism group must necessarily be of either affine or almost simple type. This result was extended by Zieschang~\cite{Zies1988} in 1988, who established that for any $2$-design satisfying $(r,\lambda)=1$, the flag-transitive automorphism group $G$ is similarly constrained to these two types. Further refinement was achieved by Zhou et al.~\cite{ZhouZhan,LZZ}, who proved that this classification holds more generally for flag-transitive automorphism groups of $2$-designs with $\lambda \geq (r,\lambda)^2$.

\textbf{Cases (b) and (c):} 
Subsequent research~\cite{ZZ2021, Zhao2022} has revealed that for $2$-designs satisfying $(r-\lambda,k)=1$, any flag-transitive automorphism group $G$ must again be of affine or almost simple type. A significant advancement was made by Zhao~\cite{Zhao2022} in 2022, who established that under the more general condition $r > \lambda(k-3)$, the flag-transitive automorphism group $G$ becomes a primitive permutation group that may additionally admit product type structure alongside the affine and almost simple types.

\textbf{Case (d):} 
The classification problem takes an interesting turn when considering designs with $(v-1,k-1) \leq 2$. As shown in~\cite[Corollary 4.6]{Kantor}, the flag-transitive automorphism group $G$ in such cases must either act $2$-transitively on points or exhibit rank $3$ with $\frac{3}{2}$-transitivity. Through application of~\cite[Theorem 1.1]{BGLPS}, these conditions necessarily imply that $G$ is of affine or almost simple type. The authors further extend this analysis to cases where $(v-1,k-1) = 3$ or $4$ by considering the broader condition $(v-1,k-1) \leq \sqrt{v-1}$. This comprehensive approach culminates in the application of reduction theorems for flag-transitive $2$-designs to the special case where $k$ is prime.

 Similar to the study of flag-transitive 2-designs,  we focus on half-flag-transitive 2-designs in this paper.

\subsection{Results}
Here, we state the results that we will prove in this paper. Firstly, we give some necessary conditions for a $2$-design to have a half-flag-transitive point-primitive automorphism group:
\begin{theorem}\label{primitive}
Let $\mathcal{D}=(\mathcal{P},\mathcal{B})$ be a $2$-$(v,k,\lambda)$ design, and  $G$ be a half-flag-transitive automorphism group of $\cal D$. Then  each of the following conditions implies that $G$ is point-primitive:
\begin{itemize}
  \item [\rm(i)] $\lambda\ge (r,2\lambda)^2$;
  \item [\rm(ii)] $r>4\lambda(k-2)$.
  \item [\rm(iii)] $(v-1,2k-2)\le2$
\end{itemize}
\end{theorem}

Next, for the case $\lambda\ge (r,2\lambda)^2$, we have the following reduction theorem for primitive groups:

\begin{theorem}\label{reduction} Let $\mathcal{D}=(\mathcal{P},\mathcal{B})$ be a $2$-$(v,k,\lambda)$ design, and $G$ be a half-flag-transitive automorphism group of $\cal D$. If $\lambda\ge (r,2\lambda)^2$, then $G$ is of  affine type, almost simple type or product type with $rank(G)=3$ and $G\le K\wr S_2$ acting on ${\cal P} = \Delta\times \Delta$ , where $K$ is a $2$-transitive group on $\Delta$ with odd size.
\end{theorem}

Finally, we analyse the case in Theorem \ref{reduction}  where $G$ is of almost simple type with sporadic socle, and we get the following:
\begin{theorem}\label{sporadic} Let $\mathcal{D}=(\mathcal{P},\mathcal{B})$ be a $2$-$(v,k,\lambda)$ design with $\lambda\ge (r,2\lambda)^2$, and $G$ be a half-flag-transitive automorphism group of $\cal D$. If $\mathrm{Soc}(G)$ is a sporadic simple group, then $G\cong {\rm HS}$ and $\cal D$ is either the unique $2$-$(176, 128, 15240)$ design or the unique $2$-$(176, 160, 19080)$ design.
\end{theorem}

\subsection{Examples}
Here we present two half-flag-transitive point-imprimitive $2$-designs with block sizes 4 and 6 respectively.
\subsubsection{2-(25,4,2) design}
Take $\mathcal{P}=\{1,2,\ldots,25\}$ to be the set of points, and $B=\{1, 2, 9, 24\}$ be a block. Now take the set of blocks to be ${\cal B}=B^G$ where $G={\tt TransitiveGroups(25,32)}$ the 32-th transitive group in the list of the {\sc Magma} \cite{Magma}-library of transitive groups with degree 25. Then $\mathcal{D}=({\cal P},{\cal B})$ is a 
$2$-$(25, 4, 2)$ design with 100 blocks. Also, the set-wise stabilizer  $G_B\cong \mathbb{Z}_2\times \mathbb{Z}_2$ has 9 orbits on the point set ${\cal P}$, two of which are $\{1,2\}$ and $\{9,24\}$. Thus, the automorphism group $G$ is half-flag-transitive and point-imprimitive.
\subsubsection{(16,6,2) biplane}
 There are exactly three non-isomorphic $(16,6,2)$  biplanes \cite{Hussain}. One arises from a difference set in \( \mathbb{Z}_2^4 \): Take the set of points \( {\cal P} = \mathbb{Z}_2^4 \), and the set of blocks \( \mathcal{B} = \{B_0 + p \mid p \in {\cal P}\} \), where 
\[
B_0 = \{\overline{0}, e_1, e_2, e_3, e_4, \sum_{i=1}^4 e_i\},
\]
and \( e_i \) is the vector with 1 in the \( i \)-th place, and 0 elsewhere, so \(\{e_1, e_2, e_3, e_4\}\) is the canonical basis for \( \mathbb{Z}_2^4 \). The automorphism group $G$ is \( 2^4 \rtimes \mathbb{Z}_3 \). Since the stabiliser \( G_{B_{0}} \cong \mathbb{Z}_3 \) has two orbits of length 3 on the six points of \( B_{0} \), and the group of translations \( 2^4 \) acts regularly on the blocks of \( {\cal B} \), \( G \) is half-flag-transitive and point-imprimitive.
\section{ Lemmas}
Here we present some auxiliary lemmas that will be needed in what follows.
\begin{lemma}\label{inq}
  Let \(\cal D\) be  a $2$-$(v,k,\lambda)$ design with \(\lambda\ge(r,2\lambda)^2 \). Then   $\frac{r^2}{(r,2\lambda)^2}>v$.
\end{lemma}
{\bf Proof.} By equations (\ref{2flagnum})  and (\ref{fisherineq}), $r^2>\lambda v\ge (r,2\lambda)^2v$ and so the result holds.  $\hfill\square$

The following lemma implies that both $r$ and $k$ are even.
\begin{lemma}  \label{eqv}
Let  $\mathcal{D}=(\mathcal{P}, \mathcal{B})$ be a $2$-design, and $\alpha \in \mathcal{P}$. If $G\le\mathrm{Aut}(\mathcal{D})$ acts half-flag-transitively on $\cal D$ then $G$ acts transitively on $\mathcal{P}$, and the point stabilizer $G_{\alpha}$ has two orbits on $P(\alpha)$ (the set of blocks through $\alpha$) with same size $\frac{r}{2}$.
\end{lemma}
\textbf{Proof.} Let $(\alpha,B)$ be a flag of $\cal D$, then
$$\frac{bk}{2}=[G:G_{\alpha B}]=[G:G_B][G_B:G_{\alpha B}]=[G:G_{\alpha}][G_{\alpha}:G_{\alpha B}].$$ 
Thus, the result can be easily derived from the Block's Lemma and Equation (\ref{flagnum}).  $\hfill\square$

\begin{lemma}\label{div} Let $\mathcal{D}=(\mathcal{P},\mathcal{B})$ be a $2$-$(v,k,\lambda)$ design, and let $G$ be a half-flag-transitive automorphism group of $\cal D$. Then $G$ is not regular, and  $\frac{r}{(r,2\lambda)}$ divides each nontrivial subdegree $d$ of $G$.
\end{lemma}
\textbf{Proof.} Suppose that $G_{\alpha}=1$ for any $\alpha\in {\cal P}$. Lemma \ref{eqv} yields that $G_{\alpha}$ has two orbits on $P(\alpha)$  with same size $\frac{r}{2}$. Thus $r=2$, a contradiction.

Let $\Gamma$ be a nontrivial $G_{\alpha}$-orbit on $\cal P$. Suppose that  $G_{\alpha}$ has two orbits $P_1$ and $P_2$ on $P(\alpha)$. Set $\mu_1=|\Gamma\cap B_1|$ and  $\mu_2=|\Gamma\cap B_2|$, where $B_i\in P_i$ ($i=1,2$). The half-flag-transitivity of $G$ and Lemma \ref{eqv} yields that $\mu_i$ is dependent of choice of $B_i$. Now we count the size of  flag set $\{(\gamma,B)\mid \alpha\in B \,{\rm and}\, \gamma\in \Gamma\}$ in two ways, and get
$$\frac{r}{2}(\mu_1+\mu_2)=\lambda|\Gamma|.$$
Therefore,   $\frac{r}{(r,2\lambda)}$ divides $|\Gamma|$.  $\hfill\square$
\begin{lemma}\label{OUT}
 There does not exist a non-abelian simple group $T$ such that $ |T| <6|\mathrm{Out}(T)|^2$.
\end{lemma}
{\bf Proof.} Assume the contrary, by \cite[Lemma 2.3]{TZ}, a non-abelian simple group $T$ satifying the inequality $ |T| <400|\mathrm{Out}(T)|^2$ then \( T \) is isomorphic to one of the following groups:
\[
PSL(2,q)\,\, \text{for }\, q=5,7,8,9,11,13,16,27, \, \,\text{or}\,\,  PSL(3,4).
\]
The facts \(|\mathrm{Out}(PSL(2,q))| = 2\) for \(q \in \{5, 7, 8, 11, 13, 16, 27\}\), \(|\mathrm{Out}(PSL(3,4))| = 12\) and \(|\mathrm{Out}(PSL(2,9))| = 4\) imply that \(|T| > 6|\mathrm{Out}(T)|^2\), a contradiction.   $\hfill\square$

\begin{lemma}\label{large}
  Let ${\cal D}$ be a  $2$-$(v,k,\lambda)$ design with a half-flag-transitive automorphism group $G$. Then
$$|G|\le4|G_{\alpha}|^3.$$ 
\end{lemma}
{\bf Proof.} The result follows from the fact that $r^2>v$ and Lemma \ref{eqv}. $\hfill\square$
\section{Primitivity}
We now prove the first result (Theorem \ref{primitive}), which is established through two propositions in this section.
\begin{proposition}\label{imprim} Let $\mathcal{D}=(\mathcal{P},\mathcal{B})$ be a $2$-$(v,k,\lambda)$ design, and $G$ be a half-flag-transitive automorphism group of $\cal D$. If  $G$ is point-imprimitive, then 
\begin{itemize}
  \item [\rm(i)] $\lambda< (r,2\lambda)^2$;
  \item [\rm(ii)] $r\le4\lambda(k-2)$.
\end{itemize}
\end{proposition}
{\bf Proof.} 
 Let $\{C_1, C_2,\ldots,C_s\}$ be an imprimitive partition system of $G$ on point set $\cal P$, where $|C_i|=w\ge2$ ($i\in\{1, 2,\ldots,s\}$) and $w\ge 2$. Then $v=sw$. Define $$\mathcal{B}_i=\{B\cap C_i\mid B\in \mathcal{B} \text{ and } |B\cap C_i|\ge1\}.$$
 The half-flag-transitivity implies that if  $\overline{B}\in \mathcal{B}_i$ then $|\overline{B}|=k_1$ or $k_2$  with $1\le k_1\le k_2\le w-1$. Thus, the substructure $\mathcal{D}_i=(\mathcal{P},\mathcal{B}_i)$ is a $2$-$(v,\{k_1,k_2\},\lambda)$ design. From equation (\ref{gl2}), we have 
 \begin{equation}\label{rlw}
 \frac{r}{2}(k_1+k_2-2)=\lambda(w-1).
 \end{equation} 
As $\lambda(s-1)=\lambda(v-1)-\lambda(v-s)$, it follows from $v=sw$ and $r(k-1)=\lambda(v-1)$ that 
\begin{equation}\label{rld}
\lambda(s-1)=\frac{r}{2}[2k-2-s(k_1+k_2-2)]=\frac{r}{2}(2s-2-D),
\end{equation}
where $D=s(k_1+k_2)-2k$.
As \(\lambda < r\) and \(\lambda (s-1)>0\) then $0<D<2s-2$.  Next, let $l=(\frac{r}{2},\lambda)$, it is clear from equation (\ref{rld}) that \( \frac{r}{2l}\mid(s-1) \).  We put \( s-1 = \frac{hr}{2l} \), where \( h \) is some positive integer. Fisher's inequality yields \( r \geq k \), hence 
\begin{align*}
l= \frac{rh}{2(s-1)}& \geq \frac{hk}{2(s-1)} = \frac{hs(k_1+k_2)}{4(s-1)} - \frac{hD}{4(s-1)}
 > \frac{h(k_1+k_2)}{4} - \frac{h}{2}\\ &= \frac{h(k_1+k_2-2)}{4} \geq \frac{k_1+k_2-2}{4},
\end{align*}
    where \( 0<D < 2s-2 \) is used. This proves that \(4l> k_1+k_2-2\). As \( \frac{r}{2l} \) divides \( w-1 \) by equation (\ref{rlw}), we have \(\frac{2(w - 1)l }{r}\geq 1\), whence \[\lambda\le\frac{2\lambda(w - 1)l }{r}\leq l(k_1+k_2-2)<4l^2=(r,2\lambda)^2.\] This proves the first claim.

   Finally, since $\cal D$ is 2-design, if $k_1=k_2$, then $k_1\ge2$ and so $k_1+k_2\ge4$; if $k_1\ne k_2$, then necessarily $1\le k_1<k_2$ and so $k_1+k_2\ge3$.  Equation (\ref{rld}) implies that \(\lambda(s-1) \geq \frac{r}{2} > 0\) and 
    \[
    2(k-1) \geq s(k_1+k_2-2) + 1 = s(k_1+k_2-3) + s + 1 \geq s - 1 + 2 \geq \frac{r}{2\lambda}+2.
    \]
Therefore, $4\lambda(k-2)\ge r$, this proves the second claim.
$\hfill\square$

\begin{proposition}\label{imprim2} Let $\mathcal{D}=(\mathcal{P},\mathcal{B})$ be a $2$-$(v,k,\lambda)$ design, and $G$ be a half-flag-transitive automorphism group of $\cal D$. If  $(v-1,2k-2)\le2$ then $G$ is point-primitive.
\end{proposition}
{\bf Proof.}\,  Equation (\ref{2flagnum}), gives 
\[
\frac{2\lambda}{(r,2\lambda)} \cdot \frac{v-1}{s} = \frac{r}{(r,2\lambda)} \cdot \frac{2k-2}{s},
\]
where $s=(v-1,2k-2)$. From the facts 
\[
\left(\frac{v-1}{s}, \frac{k-1}{s}\right) = 1 \quad \text{and} \quad \left(\frac{r}{(r,2\lambda)}, \frac{2\lambda}{(r,2\lambda)}\right) = 1,
\]
then 
\[
\frac{v-1}{s} = \frac{r}{(r,2\lambda)} \quad \text{and} \quad \frac{2k-2}{s} = \frac{2\lambda}{(r,2\lambda)}.
\]
By Lemma \ref{div}, we have 
\[
\frac{v-1}{s} \mid d. 
\] 

If $s=1$ then $d=v-1$ by Lemma \ref{div}, and so $G$ is 2-transitive also point-primitive. 

If $s=2$  then subdegrees  of $G$ is one of the following as $G$ is not regular:
\[
\left\{1, \frac{v-1}{2}, \frac{v-1}{2}\right\} \text{or} \left\{\underbrace{1,1,\cdots,1}_{\frac{v+1}{2}}, \frac{v-1}{2}\right\}. 
\]
Neither case mentioned above is possible when $G$ acts as a point-imprimitive group.$\hfill\square$
\bigskip

%
%

From now on, we assume that the following hypothesis holds: 

\textbf{Hypothesis:} Let \(\cal D\) be  a $2$-$(v,k,\lambda)$ design with \(\lambda\ge(r,2\lambda)^2 \) which admits a half-flag-transitive automorphism group \(G\).
\section{Reduction}
In this section, we establish Theorem \ref{reduction} through three key propositions. Suppose that \( (\mathcal{D},G )\) satisfy Hypothesis.  Then $G$ is point-primitive by Proposition \ref{imprim}.  According to the O'Nan-Scott theorem \cite{LPS}, primitive group $G$ is  of affine type, almost simple type, product action type, simple diagonal type, or twisted wreath product type. We adapt the method of \cite{ZhouZhan} to our setting where $G$ is half-flag-transitive. Although the flag-transitivity assumption in \cite{ZhouZhan} is relaxed to half-flag-transitivity here, the core strategy remains applicable.
%
%
%
%
\subsection{Simple diagonal type}
First suppose that \( G \) is of simple diagonal type. Let \( N =\text{Soc}(G)\cong T_1 \times \cdots \times T_m \), where each \( T_i \) (\( i \in \{1, 2, \ldots, m\} \)) is a non-abelian simple group and \( T_i \cong T \). Then $G$ is a subgroup of $W$, where
\[
W = \{(a_1, \ldots, a_m) \pi \mid a_i \in \text{Aut}(T), \pi \in S_m, \, a_i \equiv a_j\pmod{\text{Im}(T)} \text{ for all } i, j\}.
\]

\begin{proposition}\label{nsd}
Let \( (\mathcal{D},G )\) satisfy Hypothesis. Then \( G \) is not of simple diagonal type.
\end{proposition}
{\bf Proof.}
For a point \( \alpha \in {\cal P} \), we have
\[
T \cong \{(a, \ldots, a) \mid a \in T\} = N_\alpha \leq G_\alpha \leq W_\alpha,
\]
where \( W_\alpha = \{(a, \ldots, a) \pi \mid a \in \text{Aut}(T), \pi \in S_m\} \cong \text{Aut}(T) \times S_m \). Clearly,  \( v = |N : N_\alpha| = |T|^{m-1} \). 
From the proof of \cite[Proposition 3.1]{ZhouZhan}, there is a nontrivial orbit $\Delta$ of $G_{\alpha}$ on $\cal P$ such that 
$|\Delta| \le m|T|$. Thus, by Lemmas \ref{inq} and \ref{div}
\[
|T|^{m-1} = v < \frac{r^2}{(r,2\lambda)^2} < m^2 |T|^2.
\]
It follows that \( m = 2 \) or \( m = 3 \).

Since \( G_\alpha \leq \mathrm{Aut}(T) \times S_m \) and \( \frac{r}{2} \mid |G_\alpha| \), then \(\frac{ r}{2} \) divides \( |T||\mathrm{Out}(T)|m! \). From equation (\ref{2flagnum}), we have 
\[
\frac{r}{(r,2\lambda)}(2k-2) = \frac{2\lambda}{(r,2\lambda)}(|T|^{m-1}-1).
\]
It follows that $\frac{r}{(r,2\lambda)}$ divides $|T|^{m-1}-1 $, that is, \( \left(\frac{r}{(r,2\lambda)},|T|\right)=1 \), which implies $\frac{r}{(r,2\lambda)}$ divides $|\mathrm{Out}(T)|m!$. Therefore, by  Lemmas \ref{inq} and \ref{div}, we get
\[
|T|^{m-1} = v <  \frac{r^2}{(r,2\lambda)^2} < |\mathrm{Out}(T)|^2 m!^2.
\]
Thus \( |T| < 4|\mathrm{Out}(T)|^2\le  6|\mathrm{Out}(T)|^2\) if \( m=2 \); \( |T| < 6|\mathrm{Out}(T)|\le 6|\mathrm{Out}(T)|^2\) if \( m=3 \), contrary to Lemma \ref{OUT}.   $\hfill\square$

\subsection{Twisted wreath type}
Suppose that \( G \) is of twisted wreath type. Then \( N = \mathrm{Soc}(G)\cong T_1 \times \cdots \times T_m \) with \( m \geq 6 \) \cite{LPS}, where each \( T_i \cong T \) (\( i \in \{1, 2, \dots, m\} \)) is a nonabelian simple group. Also, \( N \) is regular on point set \({\cal P} \). Let \( \alpha \in {\cal P} \), then \( G = N \rtimes G_{\alpha} \)  and \( v = |T|^m \). 
\begin{proposition}\label{ntw}
Let \( (\mathcal{D},G )\) satisfy Hypothesis. Then \( G \) is not of twisted wreath type.
\end{proposition}
{\bf Proof.} Define subsets
\(
\Gamma_i = \alpha^{T_i}
\)
for all \( i \in \{1, \dots, m\} \), we have that \( \Gamma_i \cap \Gamma_j = \{\alpha\} \) for any \( i \neq j \). Then, by  the proof of \cite[Proposition 3.1]{ZhouZhan}, we have that  
 \( \bigcup_{i=1}^{m} \Gamma_i - \{\alpha\} \) is the union of some nontrivial \( G_{\alpha} \)-orbits, and 
  \[
|\bigcup_{i=1}^{m} \Gamma_i - \{\alpha\}| = m(|T| - 1).
\]
By Lemma \ref{div}, we have \( \frac{r}{(r,2\lambda)} \mid m(|T| - 1) \). Thus, from Lemma \ref{inq}
\[
|T|^m=v < \frac{r^2}{(r,2\lambda)^2} < m^2|T|^2.
\]
As \( T \) is a nonabelian simple group, therefore
\[
60^{m-2} \leq |T|^{m-2} < m^2.
\]
It follows that \( m \leq 2 \), which contradicts the fact \( m \geq 6 \). $\hfill\square$

\subsection{Product type}
Next, suppose that $G$ has a product action on $\cal P$. Then there is a group $K$ with a primitive
action (of almost simple or diagonal type) on a set $\Delta$ of size $v_0 > 5$, where
\[ {\cal P} = \Delta^m, \,\, G \le K^m \rtimes S_m=K\wr S_m,\,\, \text{and} \,\, m \geq 2. \]
\begin{proposition}\label{pd}
Let \( (\mathcal{D},G )\) satisfy Hypothesis. If $G$ is of product type then $G\le K\wr S_2$ acting on ${\cal P} = \Delta\times \Delta$ with rank \( 3 \), where $K$ is a $2$-transitive group on $\Delta$ with odd size.
\end{proposition}
{\bf proof.}\, Let \( H = K \wr S_m \), and let \( S_m \) act on \( M = \{1, 2, \ldots, m\} \). Since \( G \) is half-flag-transitive, Lemma \ref{div} implies \( [G_\alpha : G_{\alpha \beta}] =|\beta^{G_{\alpha}}|\geq \frac{r}{(r, 2\lambda)} \) for any two distinct points \( \alpha, \beta \). Therefore
\begin{equation}\label{eqpro}
  [H_\alpha : H_{\alpha \beta}] \geq [G_\alpha : G_{\alpha \beta}] \geq \frac{r}{(r, 2\lambda)} = \frac{\lambda (v - 1)}{(r, 2\lambda)(k - 1)}.
\end{equation}
Let \( \alpha = (\gamma, \gamma, \ldots, \gamma), \, \gamma \in \Gamma, \, \beta = (\delta, \gamma, \ldots, \gamma), \, \gamma \neq \delta \in \Gamma \) and let \( B \cong K^m \) be the base group of \( H \). Then 
\[
B_\alpha = K_\gamma \times \cdots \times K_\gamma, \quad B_{\alpha \beta} = K_{\gamma \delta} \times K_\gamma \times \cdots \times K_\gamma.
\]
Now \( H_\alpha = K_\gamma \wr S_m \), and \( H_{\alpha \beta} \geq K_{\gamma \delta} \times (K_\gamma \wr S_{m-1}) \). Suppose \( K \) has rank \( s \) on \( \Delta \) with \( s \geq 2 \). We can choose a \( \delta \) satisfying \( [K_\gamma : K_{\gamma \delta}] \leq \frac{v_0 - 1}{s - 1} \), so that
\[
[H_\alpha : H_{\alpha \beta}] =\frac{|H_\alpha|}{|H_{\alpha \beta}|}  \leq  \frac{|K_\gamma |^m \cdot m!}{|K_{\gamma \delta}| |K_\gamma |^{m-1} \cdot (m - 1)!} \leq m \frac{v_0 - 1}{s - 1},
\]
and hence by inequality (\ref{eqpro}),

\begin{equation}\label{inqG}
\frac{\lambda (v-1)}{(r, 2\lambda)(k-1)} \leq [G_\alpha : G_{\alpha \beta}] \leq m \frac{v_0 - 1}{s-1}.
\end{equation}
So

\begin{equation}\label{inqv0s}
 \frac{v_0^m - 1}{v_0 - 1} \leq m \frac{(r, 2\lambda)(k-1)}{\lambda(s-1)}.
\end{equation}
Now \(k-1 \leq \sqrt{(r-1)(k-1)} < \sqrt{r(k-1)}<\sqrt{\lambda v}\) and $\sqrt{\lambda}\ge (r,2\lambda)$ yield

\begin{equation}\label{ieqv}
 v_0^{m-1}< \frac{v_0^m - 1}{v_0 - 1} < m v_0^{\frac{m}{2}} \frac{(r, 2\lambda)}{\sqrt{\lambda}} \leq m v_0^{\frac{m}{2}}.
\end{equation}
Hence \(m \leq 2\), or \(m = 3\) and \(v_0 < 9\).

If \( m = 3 \),  the inequality (\ref{ieqv}) reduces to \( v_0^2 + v_0 + 1 < \frac{3v_0^{\frac{3}{2}}}{s-1} \), so that \( v_0 = 5 \) or 6 and \( s = 2 \). Now, \( S_{v_0} \wr S_3 \) has subdegrees $\{1,3(v_0-1),3(v_0-1)^2), (v_0-1)^3\}$. Since $G\le S_{v_0} \wr S_3$, $v_0\in \{5,6\}$ and Lemma \ref{div} imply $\frac{r}{(r,2\lambda)} \mid v_0-1$ which contradicts the inequality  \( \frac{r^2}{(r, 2\lambda)^2}>v=v_0^3\). Hence \( m = 2 \).

From inequality (\ref{inqv0s}), we have
\[
v_0 + 1 \leq \frac{2(r, 2\lambda)(k - 1)}{\lambda(s - 1)} < \frac{2(r, 2\lambda)v_0}{\sqrt{\lambda}(s - 1)} \leq \frac{2v_0}{s - 1}.
\]
This implies \( s = 2 \). It follows that \( K \) acts 2-transitively on \(\Delta\), and \( H = K \wr S_2 \) has rank 3 with subdegrees 1, \( 2(v_0 - 1) \), \( (v_0 - 1)^2 \).
Now \( G \leq H \), so each subdegree of \( H \) is the sum of some subdegrees of \( G \) and so 
\[
\frac{r}{(r, 2\lambda)} \mid 2(v_0 - 1).
\]
If 
\(
\frac{r}{(r, 2\lambda)} \neq 2(v_0 - 1),
\)
then 
\[
\frac{r}{(r, 2\lambda)} \leq v_0 - 1,
\]
so that
\[
 v_0^2=v<\frac{r^2}{(r, 2\lambda)^2} \leq (v_0 - 1)^2,
\]
which is a contradiction. Thus 
\(
\frac{r}{(r, \lambda)} = 2(v_0 - 1).
\)
From inequality (\ref{inqG}), we obtain that \( G \) must have a subdegree \( 2(v_0 - 1) \) and it follows that \( G \) induces a 2-transitive group \( \overline{G} \leq K \) on \(\Delta\). We conclude that \( G \) itself has rank 3 on \( \cal P \) with subdegrees $\{1,  2(v_0 - 1) , (v_0 - 1)^2 \}$. Therefore, 
\[
\frac{r}{(r, 2\lambda)} \mid (v_0 - 1)^2,
\]
i.e., 2 \( \mid v_0 - 1.\) So \( v = v_0^2 \) is an odd number. $\hfill\square$

\medskip
{\bf Proof of Theorem \ref{reduction}.} It follows immediately from Propositions \ref{nsd}-\ref{pd}. $\hfill\square$

\medskip
Through {\sc Magma} computations for parameter $v_0\in\{5,7,9,11\}$ in Proposition \ref{pd}, we have not detected any such 2-design. This empirical evidence motivates the following conjecture:
\begin{conjecture}
Let $\cal D$ be a $2$-$(v,k,\lambda)$ design with $\lambda\ge(r,2\lambda)^2$. Any half-flag-transitive automorphism group of $\cal D$ is not of product type.
\end{conjecture}

\section{Almost Simple Type with Sporadic Socle}

In this section, we always assume that Let $({\cal D}, G)$ satisfy Hypothesis. Let $N = \operatorname{Soc}(G)$, if $G$ is of almost simple type, then $N$ can be broken into four cases:
\begin{itemize}
    \item[(i)] $N$ is an alternating group $A_n$ with $n \geq 5$;
    \item[(ii)] $N$ is a sporadic simple group;
    \item[(iii)] $N$ is a Lie-Chevalley finite simple group;
    \item[(iv)] $N$ is a finite simple group of exceptional type.
\end{itemize}
We only deal with the case (ii) in this paper, that is, $N$ is one of the 26 sporadic simple groups. Since for all sporadic simple groups $N$, $|\operatorname{Out}(N)|\le2 $, then $G = N$ or $N : 2$. Recall that $G$ is point-primitive, $G_{\alpha}$ is a maximal subgroup of $G$ for $\alpha \in \mathcal{P}$ \cite[Theorem 8.2]{Wielandt}. 

The lists of maximal subgroups of $N$ and $\operatorname{Aut}(N)$ given in the Atlas \cite{Atlas} are complete except for the Monster $M$. The possible candidates for maximal subgroups $H$ of $M$, which are not listed in \cite{Atlas} are given in \cite{BW}. They are of the form $X \trianglelefteq H \leq \operatorname{Aut}(X)$, where $X$ is isomorphic to one of the following simple groups: $L_2(13)$, $U_3(4)$, $U_3(8)$, $\operatorname{Sz}(8)$. 
Obviously, $4|H|^3 \leq 4|\operatorname{Aut}(X)|^3 < |M|$, which is contrary to Lemma \ref{large}. Hence, $X$ can be ruled out.

 Let $r_{\rm max}=(v-1,|G_{\alpha}|,d_1, d_2, \ldots, d_s)$, where $d_i$ is a nontrivial subdegree of $G$. Clearly,  $\frac{r}{(r,2\lambda)} \mid  r_{\rm max}$, so $$r_{\rm max}^2\ge\frac{r^2}{(r,2\lambda)^2}>v.$$
By examining each sporadic group $G$ individually, all possible cases that satisfy this inequality are listed in Table \ref{DG}. 
\begin{table}[h]
\centering
\caption{Potential degree $v$ and automorphism group $G$}
\label{DG}
\begin{tabular}{llllllc}
\toprule
Case & \( G \)    & \( G_\alpha \)    & \(v\)    & Subdegrees & \(r_{\rm max}\) &ref.\\
\midrule
1   & \( M_{11} \)    & \( M_{10} \)    & 11    & 1, 10    & 10  &  Step \ref{s1}\\
2   &                 & \( L_2{(11)} \)   & 12    & 1, 11    & 11 &   Step \ref{s3}\\
3   &                 &  \( M_{9}:2 \)  & 55    & 1, 18, 36    & 18  &  Step \ref{s1}\\
4   & \( M_{12} \)    & \( M_{11} \)    & 12    & 1, 11    & 11   & Step \ref{s3}\\
5   & \( M_{22} \)    & \( L_3{(4)} \)    & 22    & 1, 21    & 21    &Step \ref{s3}\\
6   &                & \( A_{7} \)    & 176    & 1, 70,105    & 35    &Step \ref{s3}\\
7   & \( M_{22}:2 \)    & \( L_3{(4)}:2 \)    & 22    & 1, 21    & 21   & Step \ref{s3}\\
8   & \( M_{23} \)    & \( M_{22} \)    & 23    & 1, 22    & 22    &Step \ref{s1}\\
9   &                & \( L_3{(4)}:2_2 \)    & 253    & 1, 252    & 252    &Step \ref{s3}\\
10   &                & \( 2^4:A_{7} \)    & 253    & 1, 252    & 252    &Step \ref{s3}\\
11   & \( HS \)    & \( M_{22} \)    & 100    & 1, 99    & 99   & Step \ref{s3}\\
12  &             & \( U_3{(5)}:2 \)    & 176    & 1, 175    & 175 &  Proposition \ref{HS} \\
13   & \( HS:2 \)    & \( M_{22}:2 \)    & 100    & 1, 99    & 99 &   Step \ref{s3}\\
14   & \( M_{24} \)    & \( M_{23} \)    & 24    & 1, 23    & 23 &   Step \ref{s3}\\
15 &                 & \( M_{12}:2 \)    & 1288    & 1, 495,792    & 99 &   Step \ref{s3}\\
16   & \( Co_3 \)    & \( McL:2 \)    & 276    & 1, 275    & 275 &   Step \ref{s3}\\

\bottomrule
\end{tabular}
\end{table}

Since $\cal D$ is a 2-design with $\lambda\ge(r,2\lambda)^2$  and $G$ is a half-flag-transitive automorphism group of $\cal D$, the following Step \ref{s1} must be satisfied.
\begin{step}\label{s1} All parameters $(v,b,r,k,\lambda)$ of $\cal D$ satisfy the following conditions:
  \begin{enumerate}
    \item[\rm(i)] $2 < k < v-1$, both $k$ and $r$ are even;
    
    \item[\rm(ii)]  $\frac{r}{2}\mid |G_{\alpha}|$, $b = \frac{vr}{k}$ and $b\ge v$;
    
    \item[\rm(iii)] $\lambda = \frac{r(k-1)}{v-1}$ and $\lambda\ge(r,2\lambda)^2 $.
\end{enumerate}
\end{step}
From Step \ref{s1}, we can obtain all possible parameters of the design $\cal D$. As the half-flag-transitivity of $G$, then $[G:G_B]=b$ for any block $B\in \cal B$, and the setwise stabilizer  $G_B$ partitions the point set of $B$ into exactly two orbits of size $\frac{k}{2}$. Thus, the following Steps \ref{s2} and \ref{s3} must be satisfied. 
\begin{step}\label{s2}
There exists at least one subgroup $H$ with index \( b \) in \( G \).
\end{step}


\begin{step}\label{s3}
There exists two orbits $O_1$ and  $O_2$ of \( H\) with same size \( \frac{k}{2} \).
\end{step}

From Steps \ref{s2} and \ref{s3}, the base block $B$ of $\cal D$ may be the union of $O_1$ and $O_2$. Thus, the following Step \ref{s4} must be satisfied. 
\begin{step}\label{s4}
\( |(O_1\cup O_2)^G| = b \) and $(O_1\cup O_2)^G$ is the block set of $\cal D$.
\end{step}

All cases in Table \ref{DG} can be eliminated by Steps \ref{s1}-\ref{s3}, except for Case 12. In order to provide an intuitive understanding of the above five steps, we will now consider in detail the Case 12 in Table \ref{DG} where $G\cong{\rm HS}$ and $v=176$. With the help of the computer algebra system {\sc Magma}, we can obtain 139 $b$ under Step \ref{s1} as following:\\
\{1500, 1536, 1540, 1575, 1584, 1600, 1650, 1680, 1760, 1792, 1800, 1848, 1920, 1980, 2000, 2016, 2100, 2112, 2200, 2240, 2250, 2304,2310, 2400, 2464, 2475, 2520, 2560, 2640, 2688, 2772, 2800, 2880, 3000, 3080, 3150, 3168, 3200, 3300, 3360, 3465, 3520, 3600, 3696,3840, 3960, 4000, 4032, 4200, 4224, 4400, 4480, 4500, 4608, 4620, 4800, 4950, 5040, 5280, 5376, 5544, 5600, 5760, 6000, 6160, 6300,6336, 6400, 6600, 6720, 6930, 7040, 7200, 7392, 7680, 7920, 8000, 8064, 8400, 8800, 8960, 9000, 9240, 9600, 9900, 10080, 10560, 11088,11200, 11520, 12000, 12320, 12600, 12672, 12800, 13200, 13440, 13860, 14400, 15840, 16128, 16800, 18000, 18480, 19200, 19800, 20160,21120, 22176, 23040, 24000, 25200, 26400, 26880, 27720, 28800, 31680, 33600, 36000, 36960, 38400, 39600, 40320, 50400, 55440, 57600,63360, 67200, 72000, 79200, 80640, 84000, 100800, 110880, 115200, 126000, 161280, 201600, 252000\}.

It can be verified that 15 of which satisfy Step \ref{s2}, and all possible orbit lengths  of subgroups $H$ with index $b$ of ${\rm HS}$ are listed in Table \ref{subdegrees}:

\begin{table}[ht]
\centering
\caption{Orbit lengths of subgroups $H$ with index $b$ in the group ${\rm HS}$.}
\label{subdegrees}
\begin{tabular}{llc|llc}
\toprule
$b$  & orbit lengths & Ref. & $b$  & orbit lengths & Ref. \\ \midrule
176  & 1, 175 & Step \ref{s3} &  7700 & 80, 96 & Step \ref{s3} \\ 
  & 50, 126 & Step \ref{s3} & 8800  & 1, 7, 42, 126 & Step \ref{s3} \\ 
352  & 1, 175 & Step \ref{s3} & 11550  & 16, 160 & Step \ref{s3} \\ 
  & 50, 126 & Step \ref{s3} & 15400  & 2, 12, 72, 90 & Step \ref{s3} \\ 
1100  & 8, 168 & Step \ref{s3} & 17600  & 1, 7, 42, 126 & Step \ref{s3} \\ 
  & 56, 120 & Step \ref{s3} &   & 15, 35, 126 & Step \ref{s3} \\ 
2200  & 8, 168 & Step \ref{s3} &   & 1, 70, 105 & Step \ref{s3} \\ 
  & 56, 120 & Step \ref{s3}  & 22176  & 1, 50, 125 & Step \ref{s3} \\
3850  & 80, 96 & Step \ref{s3} & 23100  & 16, 80, 80 & $\mathcal{D}_1$ \\ 
5600  & 11, 165 & Step \ref{s3} &   & 40, 40, 96 & Step \ref{s4} \\ 
  & 66, 110 & Step \ref{s3} & 28875  & 16, 64, 96 & Step \ref{s3} \\ 
5775  & 16, 160 & Step \ref{s3} &  & 48, 64, 64 & $\mathcal{D}_2$ \\ 
\bottomrule
\end{tabular}
\end{table}

Finally, we conduct a detailed analysis of the cases in Table \ref{subdegrees} and arrive at the following conclusion:
\begin{proposition}\label{HS}
 Let \( (\mathcal{D},G )\) satisfy Hypothesis. If ${\rm Soc}(G)$ is a sporadic simple group, then $G\cong{\rm HS}$ and $\cal D$ is the  unique $2$-$(176,160, 19080)$ design with $23100$ blocks, or the  unique $2$-$(176, 128, 15240)$ design with $28875$ blocks.
\end{proposition}
{\bf Proof.}   We get $G$ by {\sc Magma} command {\tt PrimitiveGroup(176,4)}. All cases of $b$ in Table \ref{subdegrees} can be ruled out by Step \ref{s3} except cases $b=23100$ and $28875$.

For Case $b=23100$, using command {\tt Subgroups(G:OrderEqual:=n)}, where $n=\frac{|G|}{23100}$, $G$ has two conjugacy class of subgroups $H$ with index $23100$, namely, $H_1$, $H_2$.    Also, the orbits  of $H_1$ acting on $\cal P$ are $O_1,O_2,O_3$ with $|O_1|=16, |O_2|=|O_3|=80$, and  $H_2$ acting on $\cal P$ are $O_4,O_5,O_6$ with $|O_4|=|O_5|=40,|O_6|=96$. Thus, the possible values $k$ are restricted to 160 and 80. 

Suppose that $B=O_4\cup O_5$, it is clearly that $|B^G|=3850\ne 23100$ from command {\tt Orbits($B^{\wedge}G$)}. That is, this case consequently fails to satisfy Step \ref{s4}. Now assume that $B_1=O_2\cup O_3$, it is clearly that $|B_1^G|=23100$  and  $G_{B_1}$ has three orbits with size \{16,80,80\}. As the 2-transitivity of $G$, then  $\mathcal{D}_1=(\mathcal{P}, B_1^G)$ must be a half-transitive $2$-$(176,160, 19080)$ design with $23100$ blocks.

For Case $b=28875$, similarly, $G$ has two conjugacy class of subgroups $H$ with index $28875$, namely, $H_3$, $H_4$.  There is no orbits with same size when $H_3$ acting on $\cal P$, thus  $H_3$ can not be stabilizer of a block in $G$. The orbits  of $H_4$ acting on $\cal P$ are $O_7,O_8,O_9$ with $|O_7|=48, |O_8|=|O_9|=64$. Let $B_2=O_8\cup O_9$, it is clearly that $|B_2^G|=28875$  and  $G_{B_2}$ has three orbits with size \{48,64,64\}. As the 2-transitivity of $G$, then  $\mathcal{D}_2=(\mathcal{P}, B_2^G)$ must be a half-transitive $2$-$(176,128, 15240)$ design with $28875$ blocks.  $\hfill\square$

Thus, we have completed the proof of Theorem \ref{sporadic}. $\hfill\square$

\section*{Conflict of interest statement}
The authors declare no conflict of interest.

\section*{Acknowledgements}
This work is supported by the National Natural Science Foundation of China (No. 12361004) and the Natural Science Foundation of Jiangxi Province (No. 20242BAB25005).

\end{document}